\documentclass[a4paper,11pt]{article}
\usepackage{latexsym}
\usepackage{amssymb}
\usepackage{enumerate}
\usepackage[shortlabels]{enumitem}
\usepackage{xcolor}
\usepackage{commath}
\usepackage{comment}
\usepackage{amsfonts}
\usepackage{amsmath,amsthm}
\usepackage{hyperref}
\usepackage{algorithm}
\usepackage{algorithmic}
\usepackage{relsize}
\usepackage[title]{appendix}
\usepackage{framed}
\usepackage{boites}
\usepackage[pdftex]{graphicx}

\newenvironment{@abssec}[1]{%
    \if@twocolumn

      \section*{#1}%
    \else

      \vspace{.05in}\footnotesize
      \parindent .2in
 {\upshape\bfseries #1. }\ignorespaces
    \fi}

    {\if@twocolumn\else\par\vspace{.1in}\fi}

\newenvironment{keywords}{\begin{@abssec}{\keywordsname}}{\end{@abssec}}

\newcommand\keywordsname{Key words}
\newcommand\AMSname{AMS subject classifications}
\newcommand\AMname{AMS subject classification}
\newcommand\restr[2]{{% we make the whole thing an ordinary symbol
\left.\kern-\nulldelimiterspace % automatically resize the bar with \right
#1 % the function
\vphantom{|} % pretend it's a little taller at normal size
\right|_{#2} % this is the delimiter
}}
\newtheorem{theorem}{Theorem}
\newtheorem{lemma}[theorem]{Lemma}

\newcommand{\RR}{\mathbb{R}}

\def\XXint#1#2#3{{\setbox0=\hbox{$#1{#2#3}{\int}$}
\vcenter{\hbox{$#2#3$}}\kern-.5\wd0}}

\newcommand{\link}{\mathop{\circ\kern-.35em -}}

\newcommand{\ol}{\overline}
\newcommand{\pa}{\partial}

\newcommand{\dv}{\mathop{\mathrm{div}}}

\newcommand{\al}{\alpha}
\newcommand{\be}{\beta}
\newcommand{\ga}{\gamma}

\newcommand{\Si}{\Sigma}
\newcommand{\te}{\theta}

\newcommand{\Om}{\Omega}
\newcommand{\rn}{{\mathbb{R}}^N}

\newcommand\setbld[2]{\left\{ #1 \;\middle |\; #2\right\}}

\numberwithin{equation}{section}

\title{\bf Pythagorean Theorem, Law of Sines and Law of Cosines: alternative proofs via shape derivatives
}
\author{Lorenzo Cavallina 
\thanks{This research was partially supported by JSPS KAKENHI Grant Number JP23H04459, and JP22K13935, JP21KK0044.
}}
\date{}

\begin{document}

\maketitle

\begin{abstract}
We provide an alternative unified approach for proving the Pythagorean theorem (in dimension $2$ and higher), the law of sines and the law of cosines, based on the concept of shape derivative. The idea behind the proofs is very simple: we translate a triangle along a specific direction and compute the resulting change in area. Equating the change in area to zero yields the statements of the three aforementioned theorems.
\end{abstract}


%\begin{AMS}
%\end{AMS}

\pagestyle{plain}
\thispagestyle{plain}

\section{Introduction}

Let $T$ denote a triangle in the Euclidean plane with vertices $A$, $B$, and $C$. Also, let $\al$, $\be$, and $\ga$ denote the respective angles. Moreover, be a slight abuse of notation, we will use the letters $a$, $b$, and $c$ to denote both the sides $BC$, $AC$, $AB$, and their respective lengths. For instance, we will write $\int_a 1 = a$. Finally, $n_a$, $n_b$, and $n_c$ will denote the outward unit normal vectors to the sides $a$, $b$, and $c$ respectively.

In what follows, we will give an alternative proof of the following three famous results of Euclidean geometry.
\begin{enumerate}[(i)]
    \item $\ga=\pi/2\implies a^2+b^2=c^2.$ \hfill (\emph{The Pythagorean Theorem}). 
    \item $\displaystyle \frac{a}{\sin\al}=\frac{b}{\sin\be}=\frac{c}{\sin\ga}.$ \hfill (\emph{The Law of Sines}). 
    \item $a^2+b^2-2ab\cos\ga=c^2.$ \hfill (\emph{The Law of Cosines})
\end{enumerate}

The main ingredient of our proofs is the concept of shape differentiation. Consider a bounded (reference) domain $\Omega$ in $\rn$ ($N\ge2$) with Lipschitz boundary $\pa\Omega$ and outward unit normal vector $n$ (defined almost everywhere on $\pa\Omega$). Also, consider a Lipschitz continuous vector field $\xi:\rn\to\rn$ with Lipschitz continuous first derivatives and define the perturbed domain as follows:
\begin{equation*}
    \Om_t:=\setbld{x+t\xi(x)}{x\in\Om},\quad \text{for small } t\ge0.
\end{equation*}
Let $f:\rn\to\RR$ be a real-valued integrable function with integrable first derivatives and consider the map $t\mapsto \int_{\Om_t} f$. This map is differentiable at $t=0$ and its derivative is given by the following lemma.

\begin{lemma}[Hadamard's formula]\label{Hadamard}
 Under the notation above, we have 
 \begin{equation*}
  \restr{\frac{d}{dt}\left(\int_{\Om_t} f\right)}{t=0}= \int_\Om \dv(f\xi) =\int_{\pa\Om} f\xi\cdot n   
 \end{equation*}
\end{lemma}

At its core, the above result is a consequence of the change of variables formula for integrals and the divergence theorem. We refer the interested reader to \cite[Theorem 5.2.2]{HP2005} for a proof in a more general setting, where the integrand $f$ may also vary with time.

In what follows, we will only use Lemma \ref{Hadamard} in the special case when the reference domain $\Om$ is a triangle $T$, $f\equiv 1$, and $\xi$ is constant. Notice that this corresponds to considering the change in area of our triangle $T$ after a translation in the direction of $\xi$. 

\begin{figure}[h]
    \centering
    \includegraphics[width=0.5\textwidth]{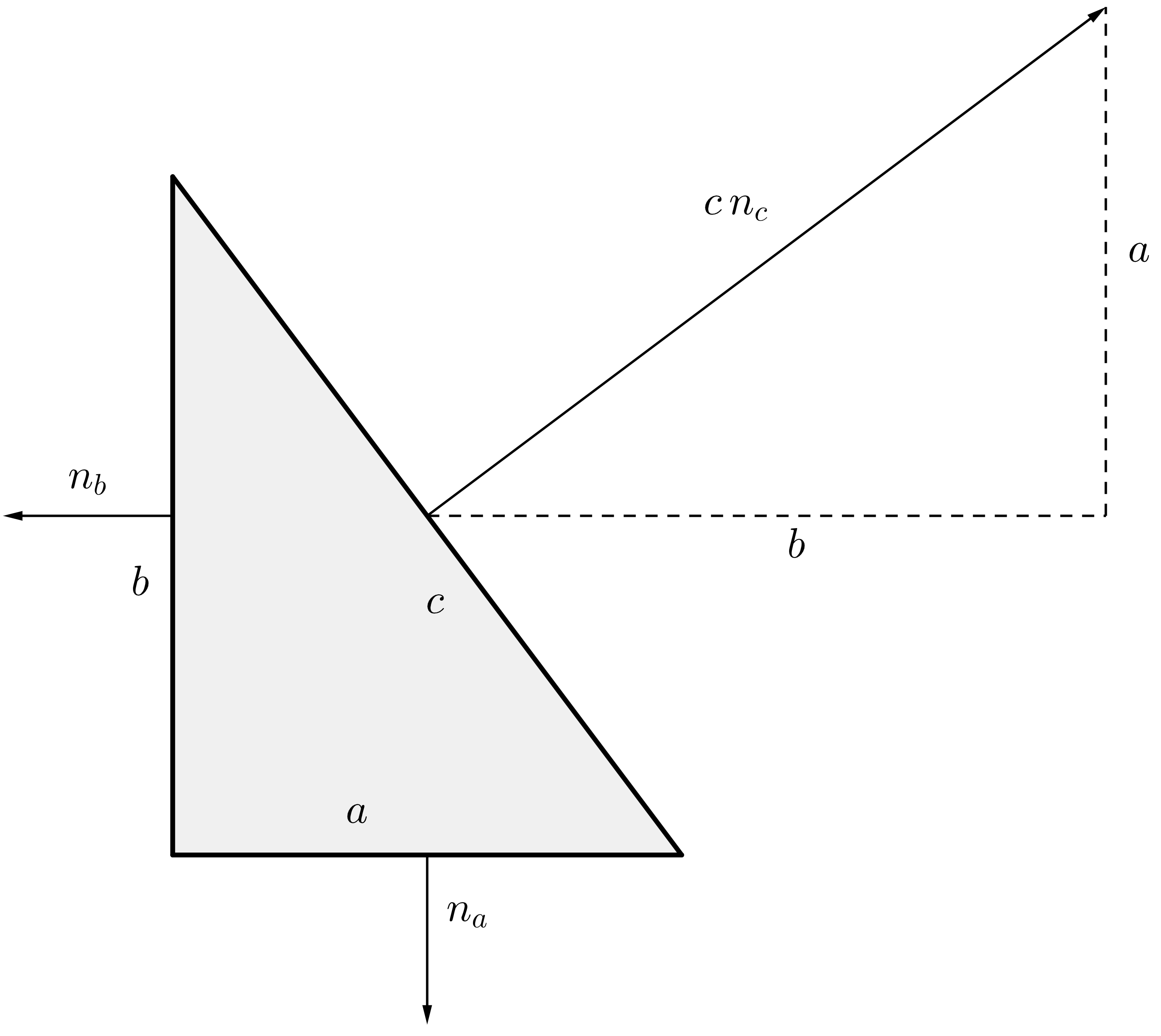}
    \caption{The construction used in our proof of the Pythagorean theorem}
    \label{fig pythagoras}
\end{figure}
 
\section{The Pythagorean theorem}

Let $T$ denote a triangle as in the introduction. Also, assume that $\gamma$ is a right angle.  
Let us consider the translation $T_t:=T+t c\, n_c$. By Lemma \ref{Hadamard} we have
\begin{equation*}
0=\restr{\frac{d}{dt}|T_t|}{t=0}=\int_{\pa T}  c\, n_c\cdot n=c\int_{c} n_c\cdot n_c + c\int_{a} n_c\cdot n_a + c\int_{b}n_c\cdot n_b=c^2-a^2-b^2.    
\end{equation*}
Here we used the fact that $c\, n_c\cdot n_a=-a$ and $c\, n_c\cdot n_b=-b$ (see Figure \ref{fig pythagoras}). 

\section{The law of sines}
\begin{figure}[h]
    \centering
    \includegraphics[width=0.4\textwidth]{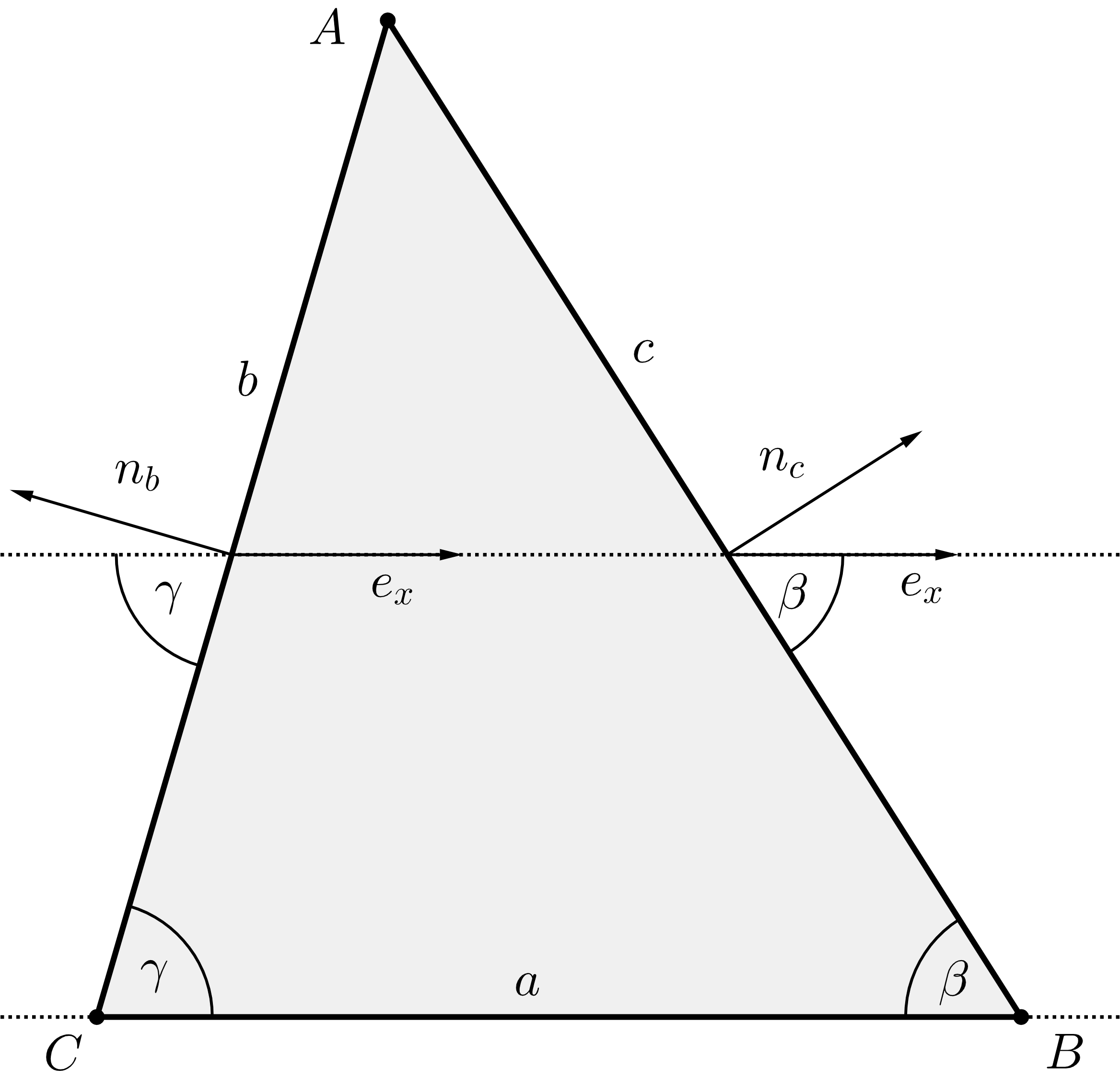}
    \caption{The construction used in our proof of the law of sines}
    \label{fig law of sines}
\end{figure}

Let $T$ denote a (non-necessarily right) triangle, as in the introduction. Also, let $e_x$ denote the unit normal vector parallel to $BC$ pointing toward $B$. Let us consider the translation $T_t:= T+te_x$. By Lemma \ref{Hadamard} we have:
\begin{equation*}
0=\restr{\frac{d}{dt}|T_t|}{t=0}=\int_{\pa T} e_x\cdot n = \int_{c} e_x \cdot n_c + \int_{b} e_x\cdot n_b = c\sin\be-b\sin\ga,
\end{equation*}
which, rearranging the terms yields the desired $b/\sin\be=c/\sin\ga$. 
Here we used the fact that $e_x\cdot n_c=\cos(\pi/2-\be)=\sin\be$ and $e_x\cdot n_b=-\cos(\pi/2-\ga)=-\sin\ga$ (see Figure \ref{fig law of sines}).
The remaining identities concerning the term $a/\sin\al$ are obtained analogously by employing a translation along the direction parallel to one of the remaining sides.  

\section{The law of cosines}
\begin{figure}[h]
    \centering    \includegraphics[width=0.6\textwidth]{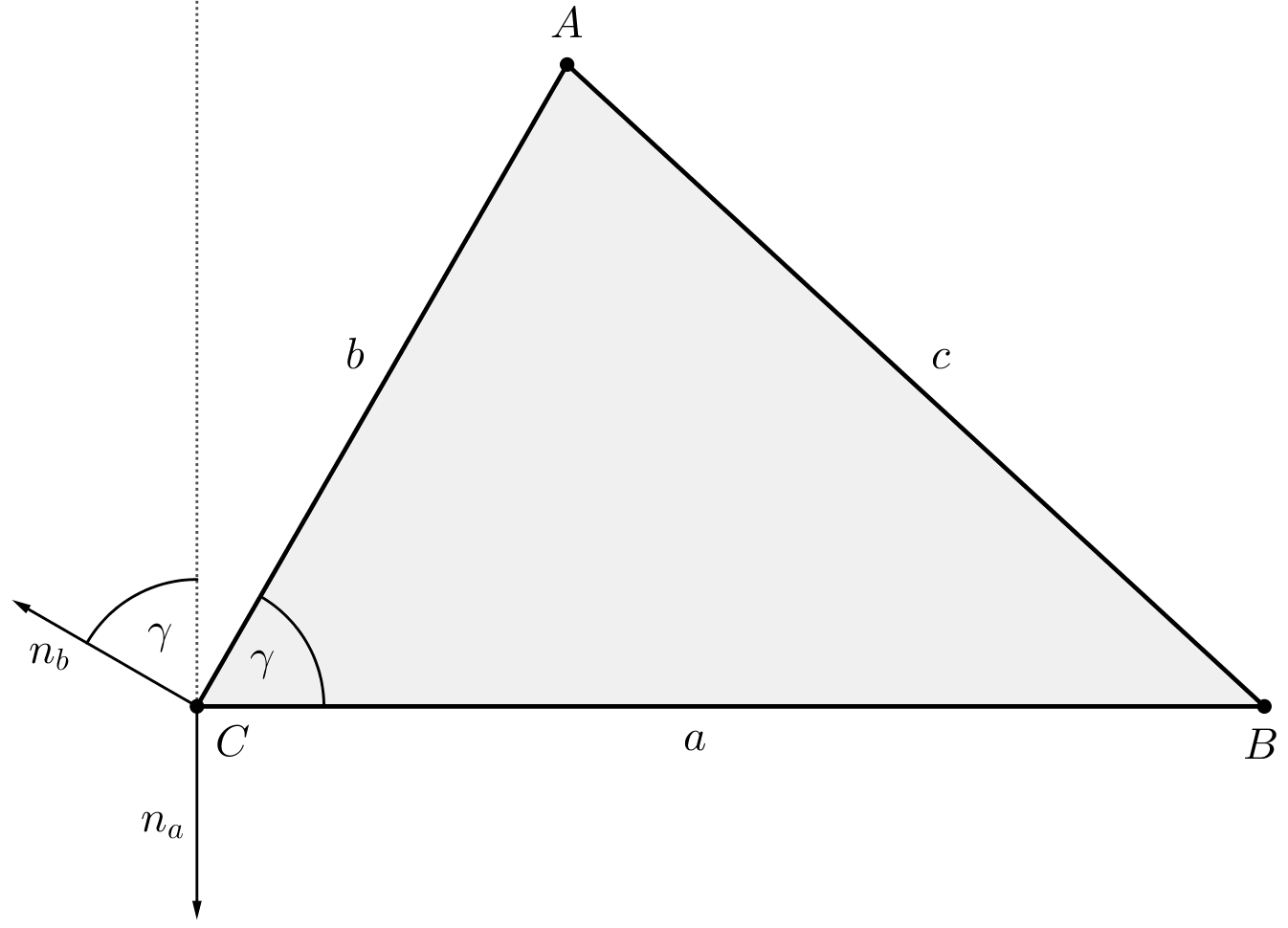}
    \caption{The construction used in our proof of the law of cosines}
    \label{fig law of cosines}
\end{figure}

Let $T$ denote a (non-necessarily right) triangle, as in the introduction. Let us consider the translation 
$T_t:=T+t(c\, n_c-a\, n_a-b\, n_b)$. By Lemma \ref{Hadamard} we have:
\begin{equation*}
    \begin{aligned}
        0=\restr{\frac{d}{dt}|T_t|}{t=0}=\int_{\pa T} (c\, n_c-a\, n_a-b\, n_b)\cdot n = \\
 c\int_{a} n_c\cdot n_a + c\int_{b} n_c\cdot n_b + c^2 -a^2 - a \int_b n_a\cdot n_b - a\int_c n_a\cdot n_c -b\int_a n_b\cdot n_a - b^2 - b\int_c n_b\cdot n_c \\
 = c^2-a^2-b^2 -2ab\, n_a\cdot n_b = c^2-a^2-b^2+2ab\cos\ga,
    \end{aligned}
\end{equation*}
which is the desired identity. Here we used the fact that $n_a\cdot n_b=-\cos\ga$ (see Figure \ref{fig law of cosines}).

\section{Concluding remarks}
Nowadays, hundreds of different proofs of the Pythagorean theorem are known (see for instance the work of Loomis \cite{Loomis}, where more than two hundred proofs are presented and analyzed). Nevertheless, it is worth mentioning that the approach given here does not seem to belong to either of the four categories introduced by Loomis.  
On a tangential note, we stress that our proofs are \emph{not trigonometric} in nature because, even though the terms $\sin\te$ and $\cos\te$ appear, they simply stand for the ratio of two lengths in a figure, and no trigonometric identity is used in a significant way (the only exception being $\cos(\pi/2-\te)=\sin\te$, which just follows from the very geometrical definitions of sine and cosine and the fact that the interior angles of a triangle in the plane add up to $\pi$ radians). 

Also, notice that we used the parallel postulate (or one  of its equivalent statements) in each proof (this cannot be avoided, as it is well-known that the parallel postulate and the Pythagorean theorem are equivalent). In particular, we explicitly made use of said postulate when chasing angles and implicitly when we used the fact that translating a pair of vectors does not change the angle between them.  

Finally, we remark that our proof of the Pythagorean theorem also generalizes to the $N$-dimensional case. The general proof in the case of right $N$-simplexes is given below (the interested reader should compare it to the one given by \cite{Eifler Rhee}, where a similar construction relying on the arbitrariness of the vector field is used instead). 
Consider an $N$-simplex $\triangle$ in $\rn$ ($N\ge2$) defined as the convex hull of $\{0,x_1,\dots, x_N\}$, where $\{x_1, \dots, x_N\}$ is an orthonormal family of vectors. For $i=1,\dots, N$, let $\Sigma_i$ denote the $(N-1)$-face containing $\{0,x_1, \dots, x_N\}\setminus\{x_i\}$, let $A_i$ denote its (($N-1$)-dimensional) area and let $n_i$ denote its outward unit normal vector. Similarly, let $\Sigma_C$ denote the remaining face, with (($N-1$)-dimensional) area $C$ and outward unit normal vector $n_C$. 

Along the same lines as section $2$, consider the translation $\triangle_t:= \triangle+tC\, n_C$. We have
\begin{equation*}
    0=\restr{\frac{d}{dt} |\triangle_t|}{t=0}=\int_{\pa\triangle} C\, n_C\cdot n = \int_{\Sigma_C} C\, n_C\cdot n_C + \sum_{i=1}^N \int_{\Sigma_i} C\, n_C\cdot n_i = C^2-\sum_{i=1}^N A_i^2,
\end{equation*}
that is the desired $N$-dimensional version of the Pythagorean theorem.
\begin{figure}[h]
    \centering    \includegraphics[width=0.6\textwidth]{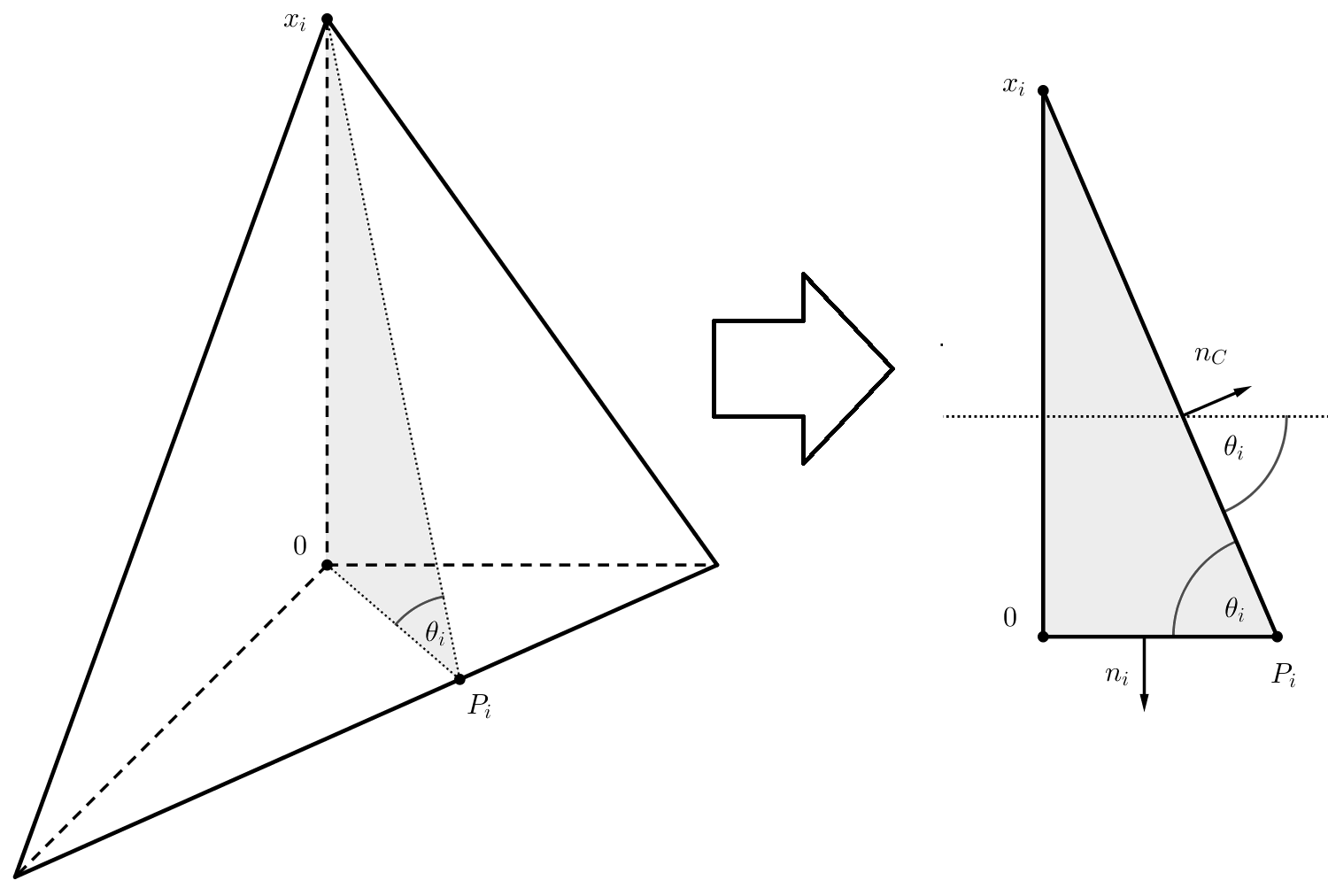}
    \caption{The construction used in the proof of the $N$-dimensional Pythagorean theorem.}
    \label{fig tridimensional}
\end{figure}
Here we used the fact that for $i=1,\dots, N$, $A_i=-C\, n_C\cdot n_i$. To show this identity, we will use the volume formula for a $k$-simplex $\triangle_k$ of base $\triangle_{k-1}$ and height $h$:
\begin{equation*}
    \left\{\text{$k$-dimensional volume of $\triangle_k$}\right\}=\frac{h}{k}\times \left\{\text{($k-1$)-dimensional volume of $\triangle_{k-1}$}\right\}.
\end{equation*}
In light of the above, since for $i=1,\dots, N$, the ($N-1$)-simplexes $\Sigma_C$ and $\Sigma_i$ share the same base $E_i$ (given by the convex hull of $\{x_1,\dots, x_N\}\setminus\{x_i\}$), it will be enough to compute the ratio of their heights. To this end, let $P_i$ be the projection of $0$ (or, equivalently, of $x_i$) onto $E_i$ and let $\te_i$ be the angle between $P_i0$ and $P_ix_i$. Since the heights of $\Si_C$ and $\Si_i$ with respect to the common base $E_i$ are given by $\ol{P_ix_i}$ and $\ol{P_i0}$ respectively, we have
\begin{equation*}
    \frac{A_i}{C}=\frac{\ol{P_i0}}{\ol{P_ix_i}}=\cos\te_i=-n_C\cdot n_i,
\end{equation*}
which is the desired identity (see also Figure \ref{fig tridimensional}).

\begin{small}

\end{small}

%\bigskip

\noindent
\textsc{
Mathematical Institute, Graduate School of Science, Tohoku University, Aoba-ku, 
Sendai 980-8578, Japan}\\
\noindent
{\em Electronic mail address:}
cavallina.lorenzo.e6@tohoku.ac.jp

\end{document}